\documentclass[12pt,a4paper]{article}
\usepackage{amsmath}
\usepackage{amscd}
\usepackage{amstext}
\usepackage{amsfonts}
\usepackage{euscript}
\usepackage{graphicx}

\def\scr{\EuScript}

\newcommand{\pcirc}{{\scriptstyle \,\circ\,}}

\newcommand{\C}{\mathbb{C}}
\newcommand{\ZZ}{\mathbb{Z}} \newcommand{\NN}{\mathbb{N}}

\DeclareMathOperator{\Hom}{\it Hom}

\DeclareMathOperator{\End}{\it End}

\DeclareMathOperator{\ann}{\rm ann}

\DeclareMathOperator{\lcm}{\rm l.c.m.}

\DeclareMathOperator{\Deg}{\rm deg}

\DeclareMathOperator{\divi}{\rm div}

\DeclareMathOperator{\rk}{\rm rk}

\DeclareMathOperator{\im}{\rm Im}

\DeclareMathOperator{\Der}{\it Der}

\DeclareMathOperator{\IC}{\rm IC}

\newcommand{\ux}{\underline{x}}

\newcommand{\uxi}{\underline{\xi}}

\DeclareMathOperator{\Sim}{{\rm Sym}}
\DeclareMathOperator{\U}{{\rm
U}}

\DeclareMathOperator{\Rees}{{\cal R}}

\DeclareMathOperator{\jac}{\rm Jac}

\newcommand{\Thetafs}{\Theta_{f,s}}

\newcommand{\bul}{\null}

\newcommand{\hol}{{\scr O}}

\newcommand{\derlogD}{\Der(\log D)}
\newcommand{\VO}{{\cal V}_0}

\newcommand{\D}{{\scr D}}
\newcommand{\K}{{\scr K}}
\newcommand{\E}{{\scr E}}
\newcommand{\F}{{\scr F}}
\newcommand{\LL}{{\scr L}}
\newcommand{\OO}{{\scr O}}

\DeclareMathOperator{\Dual}{\Bbb D}
\newcommand{\Lotimes}{\stackrel{L}{\otimes}}
\DeclareMathOperator{\DR}{DR}

\DeclareMathOperator{\Gr}{Gr}

\newcommand{\calI}{{\cal I}}

\newcommand{\OX}{{\scr O}_X}
\newcommand{\DX}{{\scr D}_X}
\newcommand{\DXp}{{\scr D}_{X,p}}
\newcommand{\OXp}{{\scr O}_{X,p}}

\DeclareMathOperator{\SP}{Sp}

\newcounter{numero}[subsection]

\newcounter{snumero}[section]

\renewcommand{\thenumero}{(\thesubsection .\arabic{numero})}
\renewcommand{\thesnumero}{(\thesection .\arabic{snumero})}
\newenvironment{corolario}{\medskip
\refstepcounter{numero}\noindent {\bf  \thenumero\ Corollary.}\
\it}{\vspace{1ex}\par}

\newenvironment{teorema}{\medskip
\refstepcounter{numero}\noindent {\bf \thenumero\ Theorem.}\
\it}{\vspace{1ex}\par}

\newenvironment{steorema}{\medskip
\refstepcounter{snumero}\noindent {\bf\thesnumero\ Theorem.}\
\it}{\vspace{1ex}\par}

\newenvironment{lema}{\medskip
\refstepcounter{numero}\noindent {\bf \thenumero\ Lemma.}\
\it}{\vspace{1ex}\par}

\newenvironment{definicion}{\medskip
\refstepcounter{numero}\noindent {\bf \thenumero\ Definition.}\
}{\vspace{1ex}\par}

\newenvironment{proposicion}{\medskip
\refstepcounter{numero}\noindent {\bf \thenumero\ Proposition.}\
\it}{\vspace{1ex}\par}

\newenvironment{nota}{\medskip
\refstepcounter{numero}\noindent {\bf \thenumero\ Remark.}\
}{\vspace{1ex}\par}

\newenvironment{ejemplo}{\medskip
\refstepcounter{numero}\noindent {\bf \thenumero\ Example.}\
}{\vspace{1ex}\par}

\newcommand\numero{\medskip\refstepcounter{numero}\noindent{\bf\thenumero}\hspace{1em}}

\newenvironment{prueba}{
\noindent {\bf Proof.}\ }{\hfill $\Box$\vspace{1ex}\par}

\title{On the logarithmic comparison theorem for integrable logarithmic connections\thanks{The authors are partially supported by MTM2004-07203-C02-01
 and FEDER.}}

\author{F. J. Calder\'{o}n Moreno and L. Narv\'{a}ez Macarro}

\date{}

\begin{document}

\maketitle

\begin{abstract}
Let $X$ be a complex analytic manifold, $D\subset X$ a free divisor
with jacobian ideal of linear type (e.g. a locally quasi-homogeneous
free divisor), $j: U=X-D \hookrightarrow X$ the corresponding open
inclusion, $\E$ an integrable logarithmic connection with respect to
$D$ and $\LL$ the local system of the horizontal sections of $\E$ on
$U$. In this paper we prove that the canonical morphisms
$$\Omega_X^{\bullet}(\log D)(\E(kD))\xrightarrow{} R j_*\LL,\quad j_!\LL
\xrightarrow{} \Omega_X^{\bullet}(\log D)(\E(-kD))$$ are
isomorphisms in the derived category of sheaves of complex vector
spaces for $k\gg 0$ (locally on $X$).
\medskip

\end{abstract}

\section*{Introduction}

\noindent Let $X$ be a $n$-dimensional complex analytic manifold. An
ideal $\calI\subset \OX$ is said to be of {\em linear type} if the
canonical homomorphism from its symmetric algebra  to its Rees
algebra is an isomorphism. We say that a divisor (=hypersurface)
$D\subset X$ is of {\em linear jacobian type} if its Jacobian ideal
is of linear type.
\medskip

This paper is devoted to prove the following result (see Corollaries
\ref{cor:main}, \ref{cor:main-bis}):\smallskip

Let $D\subset X$ be a free divisor of  linear jacobian type, $j:
U=X-D \hookrightarrow X$ the corresponding open inclusion, $\E$ an
integrable logarithmic connection with respect to $D$ and $\LL$ the
local system of the horizontal sections of $\E$ on $U$. Then, for
any point $p\in D$ there is an open neighborhood $V$ of $p$ and an
integer $k_0$ such that, for $k\geq k_0$, the restrictions to $V$ of
the canonical morphisms
$$\Omega_X^{\bullet}(\log D)(\E(kD))\xrightarrow{} R j_*\LL,\quad j_!\LL
\xrightarrow{} \Omega_X^{\bullet}(\log D)(\E(-kD))$$ are
isomorphisms in the derived category of sheaves of complex vector
spaces.
\medskip

Since any locally quasi-homogeneous free divisor is of linear
jacobian type, the above result generalizes the {\em logarithmic
comparison theorem} proved in \cite{cas_mond_nar_96} for the case
$\E=\OX$ and in  \cite{del_70}, II, \S 6, and \cite{es_vi_86},
Appendix~A, for normal crossing divisors.
\medskip

 Let us note that the Gauss-Manin construction associated with versal unfoldings
of hypersurface singularities produces non-trivial examples of
integrable logarithmic connections (with respect to the
discriminant) satisfying our hypothesis (cf.
\cite{ksaito_unifor,alek_moduli}). \medskip

Our methods are based on $D$-module theory and on our previous
results in \cite{calde_ens,calde_nar_compo,calde_nar_fourier}. See
also \cite{cas_ucha_stek,torre-45-bis,cas_ucha_exper} for related
work.
\medskip

Let us now comment on the content of this paper.
\medskip

In section 1 we introduce the notations that we will use throughout
 the paper and we recall some notions and basic results on
 Lie-Rinehart algebras, free divisors, the Bernstein construction
 and the Koszul property.  We also recall and refine
 some results in \cite{calde_nar_compo}, and focus on the
 linear type properties for a free divisor and
 the facts that any locally quasi-homogeneous free divisors is of linear jacobian type
 and that any free divisor of linear jacobian type is Koszul free.
 \medskip

 In section 2 we give an improved version of our
 characterization theorem in \cite{calde_nar_fourier} of
 the logarithmic comparison problem for
 logarithmic integrable connections. By the way, we deduce a new and short proof of the
 logarithmic comparison theorem for the trivial connection in \cite{cas_mond_nar_96}.
 \medskip

 In section 3 we state and prove the main results of this paper.
 Namely, a ``parametric" comparison theorem between logarithmic and
 usual Bernstein-Kashiwara modules associated with integrable
 logarithmic connections  (see Theorem \ref{teo:main-0} and Corollary
 \ref{cor:main-0}), and the logarithmic comparison theorem for
 integrable logarithmic connections, both with respect to free divisors of
  linear jacobian type.
 \medskip

 In section 4 we apply the above results to describe
 algebraically ``intersection $D$-modules". Namely, given an
 integrable logarithmic connection $\E$ with respect to a free divisor of
 linear jacobian type $D\subset X$, we describe in terms of $\E$ the regular
 holonomic $D$-module which corresponds, via the Riemann-Hilbert
 correspondence of Mebkhout-Kashiwara, to the intersection complex
 of Deligne-Goresky-MacPherson associated with the local system of
 horizontal sections of $\E$ on $X-D$.
 \medskip

 We thank T. Torrelli for useful comments on a previous version on
 this paper.

\section{Notations and preliminary results}

Let $X$ be a $n$-dimensional complex analytic manifold and $D\subset
X$ a hypersurface (= divisor), and let us denote by $j: U=X-D
\hookrightarrow X$ the corresponding open inclusion. We denote by
$\pi:  T^{\ast}X \to X$ the cotangent bundle, $\OX$ the sheaf of
holomorphic functions on $X$, $\DX$ the sheaf of linear differential
operators on $X$ (with holomorphic coefficients), $\Gr \DX$ the
graded ring associated with the filtration $F^{\bul}$ by the order
and $\sigma(P)$ the principal symbol of a differential operator $P$.
If $J\subset \DX$ is a left ideal, we denote by $\sigma (J)$ the
corresponding graded ideal of $\Gr \DX$.

Let us denote by $\jac(D)\subset \OX$ the Jacobian ideal of
$D\subset X$, i.e. the coherent ideal  of $\OX$ whose stalk at any
$p\in X$ is the ideal  generated by $h,\frac{\partial h}{\partial
x_1},\dots,\frac{\partial h}{\partial x_n}$, where $h\in {\scr
O}_{X,p}$ is any reduced local equation of $D$ at $p$ and
$x_1,\dots,x_n\in{\scr O}_{X,p}$ is a system of local coordinates
centered at $p$.

We say that $D$ is quasi-homogeneous at $p\in D$ if there is a system
of local coordinates  $\ux$ centered at $p$ such that the germ
$(D,p)$ has a reduced weighted homogeneous defining equation (with
strictly positive weights) with respect to $\ux$. We say that $D$ is
locally quasi-homogeneous if it is so at each point $p\in D$.

For any bounded complex $\K$ of sheaves of $\C$-vector spaces on
$X$, let us denote by $\K^{\vee}= R \Hom_{\C_X}(\K,\C_X)$ its
Verdier dual.

If $A$ is a commutative ring and $M$ an $A$-module, we will denote
by $\Sim_A(M)$ its symmetric algebra. If $I\subset A$ is an ideal,
we will denote by $\Rees (I)= \oplus_{d=0}^\infty I^d t^d \subset
A[t]$ its Rees algebra.

\subsection{Lie-Rinehart algebras}

Let $k \to A$ be a homomorphism of commutative rings.

Let us denote by $\Der_k(A)$ the $A$-module of $k$-linear
derivations $\delta:A\to A$. It is a left sub-$A$-module of
$\End_k(A)$ closed by the bracket $[-,-]$. If $\delta,\delta'\in
\Der_k(A)$ and $a\in A$ we have $ [\delta,a\delta'] =
a[\delta,\delta'] + \delta(a)\delta'$.

\begin{definicion} (Cf. \cite{MR0055323,MR0125867,rine-63})
A Lie-Rinehart algebra over $(k,A)$, or a $(k,A)$-Lie algebra, is an
$A$-module $L$ endowed with a $k$-Lie algebra structure and an
$A$-linear map $\rho: L \to \Der_k(A)$, called ``anchor map", which
is also a morphism of Lie algebras and satisfies
$$ [\lambda,a\lambda'] = a[\lambda,\lambda'] + \rho(\lambda)(a)\lambda'$$
for $\lambda,\lambda'\in L$ and $a\in A$.
\end{definicion}

 In order to simplify, we write
$\lambda(a) \stackrel{\text{not.}}{=} \rho(\lambda)(a)$ for
$\lambda\in L$ and $a\in A$.

\begin{ejemplo}  1) The first example of Lie-Rinehart algebra is $L =
\Der_k(A)$ with the identity as anchor morphism.\\
2) More generally, for any ideal $I\subset A$, the set $$
\Der_k(\log I) := \{\delta\in\Der_k(A)\ |\ \delta(I)\subset I\}$$ is
a sub-$A$-module and a sub-$k$-Lie algebra of $\Der_k(A)$ which
becomes a Lie-Rinehart algebra by considering the inclusion
$\Der_k(\log I) \hookrightarrow \Der_k(A)$ as anchor map.
\end{ejemplo}

\begin{definicion} Let $L,L'$ be Lie-Rinehart algebras over $(k,A)$.
A morphism of Lie-Rinehart algebras from $L$ to $L'$ is an
$A$-linear map $F:L\to L'$ which is a morphism of Lie algebras and
satisfies $ \lambda(a) = F(\lambda)(a)$, $\forall \lambda\in L$,
$\forall a\in A$.
\end{definicion}

\begin{definicion} An $A$-ring is a (not necessarily commutative)
ring $B$ with a ring homomorphism $\eta:A\to B$. We say that the
$A$-ring $(B,\eta)$ is central over $k$ if $\eta(c)b = b\eta(c)$ for
any $b\in B$ and any $c\in k$.
\end{definicion}

\begin{definicion} Let $L$ be a Lie-Rinehart algebra over $(k,A)$ and
$R$ a $A$-ring central over $k$. We say that a $k$-linear map
$\varphi:L\to R$ is {\em admissible} if:
\begin{enumerate}
\item[a)] $\varphi (a\lambda) = a\varphi(\lambda)$ for $\lambda\in L$
and $a\in A$, i.e. $\varphi$ is a morphism of left $A$-modules,
\item[b)] $\varphi([\lambda,\lambda']) =
[\varphi(\lambda),\varphi(\lambda')]$ for $\lambda,\lambda'\in L$,
i.e. $\varphi$ is a morphism of Lie algebras,
\item[c)] $\varphi(\lambda)a-a\varphi(\lambda) = \lambda(a) 1_R$ for
$\lambda\in L$ and $a\in A$.
\end{enumerate}
\end{definicion}

\begin{teorema} (\cite{rine-63}) For any Lie-Rinehart algebra $L$ over $(k,A)$ there
exists an $A$-ring $U$, central over $k$, and an admissible map
$\theta: L \to U$ which are universal in the sense that, for any
$A$-ring $R$ central over $k$ and any admissible map $\varphi:L\to
R$, there exists a unique $A$-ring homomorphism $h:U\to R$ such that
$h \circ \theta = \varphi$.
\end{teorema}

The pair $(U,\theta)$ in the above theorem is clearly unique, up to
a unique isomorphism. It is called the {\it enveloping algebra} of
$L$ and it is denoted by  $\U(L)$. Some authors call $\U(L)$ the
{\it universal algebra} of $L$ (cf. \cite{MR1058984}).
\medskip

The algebra $\U(L)$ has a natural filtration $F^{\bullet}$ given by
the powers of the image of $\theta$. If $L$ is a projective
$A$-module, the Poincar\'{e}-Birkhoff-Witt theorem \cite{rine-63}
asserts that its associated graded ring is canonically isomorphic to
the symmetric algebra of the $A$-module $L$, and so the map $\theta$
is injective.
\medskip

For any (commutative) scalar extension $k\to k'$ and any Lie-Rinehart algebra $L$ over $(k,A)$, $k'\otimes_k L$ inherits
an obvious Lie-Rinehart algebra structure over $(k',k'\otimes_k A)$.
\medskip

\numero \label{nume:carac-Umod} From the universal property of
$\U(L)$, any left $\U(L)$-module $M$ is determined by the
admissible map
\begin{equation} \label{eq:cone-1}
\lambda\in L \mapsto [m\mapsto \lambda m]\in \End_k(M).\end{equation}

Let us suppose for now that $L$ is a projective $A$-module of
finite rank, and let us consider the dual $A$-module $\Omega_L:=
\Hom_A(L,A)$ and the ``exterior differential"
$$d:A\to \Omega_L,\quad (da)(\lambda)= \lambda(a),\quad a\in
A.$$ The map (\ref{eq:cone-1}) is so uniquely determined by the
connection
\begin{equation*}
\nabla: M \to \Omega_L\otimes_A M,\quad \nabla(m)(\lambda) = \lambda
m,\quad m\in M, \lambda\in L,\end{equation*} where we have
identified $\Omega_L\otimes_A M = \Hom_A(L,M)$. The connection
$\nabla$ satisfies the Leibniz rule $\nabla(am) = a\nabla(m) +
(da)\otimes m$ and the admissibility of the map (\ref{eq:cone-1}) is
equivalent to the integrability condition on $\nabla$ in the usual
sense cf. \cite{del_70}, I,~2.14.
\medskip

\numero {\em The Cartan-Eilenberg-Chevalley-Rinehart-Spencer
complexes} (cf.
\cite{cheva-eilen,car_eil,rine-63,quillen-phd,mal-spencer,kas_master_trad})
\label{nume:CECRSP}

\medskip

In the following, let us suppose that $L\subset L'$ is a pair of
Lie-Rinehart algebras over $(k,A)$ and $E$ is a left $\U(L)$-module.

The Cartan-Eilenberg-Chevalley-Rinehart-Spencer complex associated
with $(L,L',E)$ is the complex $\SP_{L,L'}(E)$ defined by
$\SP^{-r}_{L,L'}(E) = \U(L') \otimes_A \left(\bigwedge^r
L\right)\otimes_A E$, $r\geq 0$ and the differential
$\varepsilon^{-r}:\SP^{-r}_{L,L'}(E)\xrightarrow{}
\SP^{-(r-1)}_{L,L'}(E)$ is given by:
\begin{eqnarray*}&
\varepsilon^{{-r}}(P\otimes(\lambda_1\wedge\cdots\wedge\lambda_r)\otimes
e) =&\\ &={\displaystyle\sum_{i=1}^r} (-1)^{i-1}
(P\lambda_i)\otimes(\lambda_1\wedge\cdots\wedge\widehat{\lambda_i}
\wedge\cdots\wedge\lambda_r)\otimes e \quad-&\\
& - {\displaystyle\sum_{i=1}^r} (-1)^{i-1}
P\otimes(\lambda_1\wedge\cdots\wedge\widehat{\lambda_i}
\wedge\cdots\wedge\lambda_r)\otimes (\lambda_i e)\quad +&\\ &+
{\displaystyle\sum_{1\leq i<j\leq r} }
(-1)^{i+j}P\otimes([\lambda_i,\lambda_j]\wedge\lambda_1\wedge\cdots\wedge
\widehat{\lambda_i}\wedge\cdots\wedge\widehat{\lambda_j}
       \wedge\cdots\wedge\lambda_r)\otimes e&
\end{eqnarray*}
for $r\geq 2$, and $\varepsilon^{{-1}} (P\otimes\lambda_1\otimes e
)=(P\lambda_1)\otimes e - P\otimes (\lambda_1 e)$ for $r=1$, and
$P\in \U(L')$, $\lambda_i\in L$, $e\in E$.
\smallskip

\noindent We also have an obvious natural augmentation
\begin{equation}
\label{eq:nat-aug} \varepsilon^0:\SP^{0}_{L,L'}(E) = \U(L')\otimes_A
E \to  h^0\left(\SP_{L,L'}(E)\right) = \U(L')\otimes_{\U(L)}
E.\end{equation}

We write $\SP_{L,L'} = \SP_{L,L'}(A)$ and $\SP_L(E) = \SP_{L,L}(E)$.
Let us note that $\SP_{L,L} = \SP_L(A)$ and
\begin{equation} \label{eq:notice} \SP_{L,L''}(E) =
\U(L'')\otimes_{\U(L')} \SP_{L,L'}(E)
\end{equation}
 for $L''\supset L'$ a third Lie-Rinehart algebra over $(k,A)$.

\begin{proposicion} \label{prop:sp-res} Let us suppose that $L$ is a projective
$A$-module of finite rank and that $E$ is a left $\U(L)$-module,
flat (resp. projective) over $A$. Then the complex $\SP_L(E)$ is a
finite $\U(L)$-resolution (resp. a finite projective
$\U(L)$-resolution) of $E$. Moreover, if $L$ and $E$ are free
$A$-modules, then $\SP_L(E)$ is a finite free $\U(L)$-resolution of
$E$.
\end{proposicion}

\begin{prueba} We proceed as in \cite{calde_nar_fourier}, p.~52. We
consider the filtration on the augmented complex $\SP_L(E) \to E$
given by
$$F^i \SP^{-k}_L(E) = (F^{i-k}\U(L)) \otimes_A
\stackrel{k}{\bigwedge} L \otimes_A E,\quad F^i E = E, \quad i\geq
0.$$ Its graded complex is, by using the Poincar\'{e}-Birkhoff-Witt
theorem, canonically isomorphic to the tensor product by $-\otimes_A
E$ of the augmented complex
\begin{equation}\label{eq:sim-aug} \Sim_A (L) \otimes_A \bigwedge^{\bullet} L
 \xrightarrow{d^0} A,
\end{equation}
where the differential is given by
$$ d^{-k}
(P\otimes(\lambda_1\wedge\cdots\wedge\lambda_k)) =\sum_{i=1}^k
(-1)^{i-1}
(P\lambda_i)\otimes(\lambda_1\wedge\cdots\wedge\widehat{\lambda_i}
\wedge\cdots\wedge\lambda_k)$$ for $P\in \Sim_A (L)$,
$\lambda_1,\dots,\lambda_k\in L$ and $k=1,\dots, \rk_A L$, and
$$d^0: \Sim_A(L) \otimes_A \stackrel{0}{\bigwedge} L =\Sim_A(L)
\xrightarrow{} A$$ is the obvious augmentation. The proposition
follows from the exactness of (\ref{eq:sim-aug}) (cf.
\cite{bou_a_10}, {\S}9,~3) and the flatness of $E$.
\end{prueba}
\medskip

\numero {\em Lie algebroids}
\medskip

The notions and results above can be easily generalized to the case
in which our ring homomorphism $k\to A$ is replaced by a morphism of
sheaves of commutative rings ${\scr K} \to {\scr A}$ on a topological
space, for instance when $X$ is a complex analytic manifold and we
consider the morphism $\C_X \to \OX$ or $\C_X[s] \to \OX[s]$. In that
case it is customary to talk about {\it Lie algebroids} instead of
Lie-Rinehart algebras. If $\scr L$ is a Lie algebroid over $({\scr
K},{\scr A})$, its stalk ${\scr L}_p$ at a point $p$ is a
Lie-Rienhart algebra over $({\scr K}_p,{\scr A}_p)$. We leave the
reader to decide the details (see \cite{macken-87,macken-05} as
general references for Lie algebroids on differentiable manifolds and
\cite{MR2000m:32015} for the interplay between complex Lie algebroids
and $\D$-module theory).

\begin{ejemplo} 1) The first example of Lie algebroid is $\LL =
\Der_{\C}(\OX)$ with the identity as anchor morphism.\\
2) The sheaf of differential operators of order $\leq 1$, $F^1 \DX =
\OX \oplus \Der_{\C}(\OX)$, with the projection $F^1 \DX \to
\Der_{\C}(\OX)$ as anchor morphism, is a Lie algebroid.\\
3) Any submodule $\LL \subset \Der_{\C}(\OX)$ which is closed for
the bracket is a Lie algebroid with the inclusion as anchor
morphism. This applies in particular to $\LL = \Der(\log D)=
\{$logarithmic vector fields with respect to $D\}$
\cite{ksaito_log}.
\end{ejemplo}

\subsection{Logarithmic derivations and free divisors}
\label{subsec:log-der}

 We say that $D$ is a {\em free divisor} \cite{ksaito_log} if the
$\OX$-module $\Der(\log D)$ of logarithmic vector fields with respect
to $D$ is locally free (of rank $n$), or equivalently if the
$\OX$-module $\Omega^1_X(\log D)$ of logarithmic 1-forms with respect
to $D$ is locally free (of rank $n$).

Normal crossing divisors, plane curves, free hyperplane arrangements
(e.g. the union of reflecting hyperplanes of a complex reflection
group), discriminant of left-right stable mappings or bifurcation
sets are example of free divisors.

Let us denote by $\DX(\log D)$ the $0$-term of the
Malgrange-Kashiwara filtration with respect to $D$ on the sheaf
$\DX$ of linear differential operators on $X$ (cf. \cite{mai-meb},
Def 4.1-1). When $D$ is a free divisor, the first author has proved
in \cite{calde_ens} that $\DX(\log D)$ is the universal enveloping
algebra of the Lie algebroid $\Der(\log D)$, and so it is coherent
and it has noetherian stalks of finite global homological dimension.
Locally, if $\{\delta_1,\dots,\delta_n\}$ is a local basis of the
logarithmic vector fields on an open set $V$, any differential
operator in $\Gamma(V,\DX(\log D))$ can be written in a unique way
as a finite sum \begin{equation} \label{eq:estruc-V0}
\sum_{\substack{\alpha \in\NN^n\\
|\alpha|\leq d}} a_{\alpha} \delta_1^{\alpha_1}\cdots
\delta_n^{\alpha_n},\end{equation}where the $a_{\alpha}$ are
holomorphic functions on $V$.

\subsection{The ring $\D[s]$ and the Bernstein construction} \label{subsec:D-bernstein}

Let $p$ be a point in $D$ and $f\in\hol=\hol_{X,p}$ a reduced local
equation of $D$. Let us write $\D= \D_{X,p}$.

On the polynomial ring $\D[s]$, with $s$ a central variable, there
are two natural filtrations: the filtration induced by the order
filtration on $\D$, that we also denote by $F^{\bul}$, and the {\em
total order} filtration $F_T^{\bul}$ given by
$$ F_T^k \D[s] = \sum_{i=0}^k (F^i \D) s^{k-i},\quad \forall k\geq 0.$$
For each $P\in \D[s]$ let us denote by $\sigma(P)$ (resp.
$\sigma_T(P)$) its principal symbol in $\Gr_{F^{\bul}} \D [s]$
(resp. in $\Gr_{F_T^{\bul}} \D [s]$).
\smallskip

The filtered ring $(\D[s],F^{\bul})$ is the ring of $\C[s]$-linear
differential operators with coefficient in $\hol[s]$ and so it is
the enveloping algebra of the Lie-Rinehart algebra
$\Der_{\C[s]}(\hol[s])=\C[s]\otimes_{\C} \Der_{\C}(\hol)$ over
$(\C[s],\hol[s])$, whereas the filtered ring $(\D[s],F_T^{\bul})$ is
the enveloping algebra of the Lie-Rinehart algebra $F^1\D = \hol
\oplus \Der_{\C}(\hol)$ over $(\C,\hol)$ whose anchor map is the
projection $\hol \oplus \Der_{\C}(\hol) \to \Der_{\C}(\hol)$. In the
latter case the canonical map $ F^1\D \to \D[s]$ sends every
$a\in\hol$ to $as$.

The canonical maps $$\eta:\Sim_{\hol[s]} (\Der_{\C[s]}(\hol[s])) \to
\Gr_{F^{\bul}} \D [s],\quad \eta_T:\Sim_{\hol} (F^1\D) \to
\Gr_{F_T^{\bul}} \D [s]$$ are isomorphisms of graded
$\hol[s]$-algebras and $\hol$-algebras respectively.

The free module of rank one over the ring $\hol[f^{-1},s]$ generated
by the symbol $f^s$, $\hol[f^{-1},s]f^s$, has a natural left module
structure over the ring $\D[s]$: the action of a derivation
$\delta\in\Der_{\C}(\hol)$ is given by $\delta (f^s) = \delta (f) s
f^{-1} f^s$ (see \cite{bern_72}).

Let us call $\jac(f)=\jac(D)_p$ the Jacobian ideal of $f$,
$$\varphi_0:\Sim_{\hol} (F^1\D) \to \Rees (\jac(f))\subset \hol[t]$$ the composition of
the canonical surjective map $\Sim_{\hol} \jac(f) \to \Rees
(\jac(f))$ with the surjective map $\Sim_{\hol} (F^1\D) \to
\Sim_{\hol} (\jac(f)) $ induced by
$$ P \in F^1\D \mapsto P(f)\in \jac(f),$$and
\begin{equation} \label{eq:varphi} \varphi := \varphi_0 \pcirc \eta_T^{-1}: \Gr_{F_T^{\bul}} \D
[s] \to \Rees (\jac(f)).\end{equation}

For each $P\in\D[s]$ of total order $d$, we have that $ P(f^s)=Q(s)
f^{-d} f^s$ where $Q(s)$ is a polynomial of degree $d$ in $s$ with
coefficients in $\hol$. Let us call $C_{P,d}\in \hol$ the highest
coefficient of $Q(s)$.

The following lemma is well-known and the proof is straightforward
(cf. \cite{yano_78}, chap.~I, Prop.~2.3).

\begin{lema} \label{lema:coro}
With the above notations,  we have $\varphi (\sigma_T(P))= C_{P,d}
t^d$ and so $ \sigma_T(\ann_{\D[s]} f^s) \subset \ker \varphi$.
\end{lema}

It is clear that $F^0_T \ann_{\D[s]} f^s = 0$ and that
\begin{equation} \label{eq:thetafs}
\Thetafs:= F^1_T \ann_{\D[s]} f^s
\end{equation}
 is formed by the
operators $\delta - \alpha s$ with $\delta\in \Der_{\C}(\hol)$,
$\alpha\in \hol$ and $\delta(f)=\alpha f$. One easily sees that the
$\hol$-linear map
$$\textstyle \delta \in \Der(\log D)_p \mapsto \delta - \frac{\delta(f)}{f}
s\in \Thetafs$$ is an isomorphism of Lie-Rinehart algebras over
$(\C,\hol)$. We obtain a canonical isomorphism
$$ \Thetafs \simeq \Gr^1_{F_T} \ann_{\D[s]} f^s.$$

On the other hand, the homogeneous part of degree one $[\ker
\varphi]_1\subset\ker \varphi$ is also canonically isomorphic to
$\Thetafs$, and so we obtain
\begin{equation*}
\Gr^1_{F_T} \ann_{\D[s]} f^s
\left( = \left[ \sigma_T(\ann_{\D[s]} f^s)\right]_1=\sigma_T(\Thetafs) \right) = [\ker
\varphi]_1.
\end{equation*}

\subsection{Divisors of linear type}

\begin{definicion}
(Cf. \cite{vascon_cmcaag}, {\S}7.2) Let $A$ be a commutative ring
and $I\subset A$ an ideal. We say that $I$ is of {\em linear type}
if the canonical (surjective) map of graded $A$-algebras $ \Sim_A(I)
\to \Rees(I)$ is an isomorphism.
\end{definicion}
\medskip

Ideals generated by a regular sequence are the first example of
ideals of linear type.

\begin{definicion} (see also \cite{simis-alg-free})
 We say that the divisor $D$ is of {\em
linear jacobian type} at $p\in D$ if the stalk at $p$ of its
jacobian ideal is of linear type. We say that $D$ is of {\em linear
jacobian type} if it is so at any $p\in D$.
\end{definicion}

\begin{nota} \label{nota:expli-clt}
To say that $D$ is of linear jacobian type at $p$ is equivalent to
saying that $\ker \varphi$ (see (\ref{eq:varphi})) is generated by
its homogeneous part of degree 1, $[\ker \varphi]_1 =
\sigma_T(\Thetafs)$.
\end{nota}
\smallskip

\noindent Theorem 5.6 of \cite{calde_nar_compo} can be rephrased in
the following way:

\begin{teorema} \label{teo:lqh->clt} Any locally quasi-homogeneous free divisor is of linear jacobian type.
\end{teorema}

\begin{definicion} Let $p\in D$ and let us write $\hol=\OXp$ and $\D=\DXp$.
We say that $D$ is of {\em differential linear type} at $p\in D$ if
for some (or any, one easily sees that this condition does not
depend on the choice of the local equation) reduced local equation
$f\in\hol$ of $D$ at $p$, the ideal $\ann_{\D[s]} f^s$ is generated
by total order one operators, i.e. (see (\ref{eq:thetafs})) $
\ann_{\D[s]} f^s = \D[s]\cdot \Thetafs$. We say that $D$ is of {\em
differential linear type} if it is so at any $p\in D$.
\end{definicion}

It is clear that the set of points at which a divisor $D$ is of
linear jacobian or differential linear type is open in $D$.

\begin{proposicion} \label{prop:clt->dlt} If the divisor $D$ is of
linear jacobian type (at $p\in D$), then it is of  differential
linear type (at $p\in D$) and if $f\in \OXp$ is a reduced local
equation of $D$ at $p$, then
$$ \Gr_{F_T} \ann_{\D[s]} f^s \left( = \sigma_T(\ann_{\D[s]}
f^s)\right) = \ker \varphi.$$
\end{proposicion}

\begin{prueba}  It is the same proof as
those of Proposition 3.2 in \cite{calde_nar_compo}, but here we
consider $ \Gr_{F_T} \ann_{\D[s]} f^s$ and the ``true'' jacobian
ideal $\jac(f)= (f,f'_{x_1},\dots,f'_{x_n})$ instead of $ \Gr_{F}
\ann_{\D} f^s$ and $J_f=(f'_{x_1},\dots,f'_{x_n})$.
\end{prueba}

\subsection{The Koszul property}

In this section, we fix a homomorphism of commutative rings $k \to
A$ (resp. a homomorphism of sheaves of commutative rings ${\scr K}
\to {\scr A}$ on a topological space $M$) and all
$(k,A)$-Lie-Rinehart algebras (resp. all Lie algebroids over $({\scr
K},{\scr A})$) will be free $A$-modules of finite rank (resp.
locally free of finite rank over $\scr A$).

We also assume that $D\subset X$ is a free divisor.
\smallskip

Let us recall that $D$ is a {\em Koszul free} divisor
\cite{calde_ens} at a point $p\in D$ if the symbols of any (or some)
local basis $\{\delta_1,\dots,\delta_n\}$ of $\Der(\log D)_p$ form a
regular sequence in $\Gr {\scr D}_{X,p}$. We say that $D$ is a {\em
Koszul free} divisor if it is so at any point $p\in D$.  Actually, as
M. Schulze pointed out, Koszul freeness is equivalent to holonomicity
in the sense of \cite{ksaito_log}.

Plane curves and locally quasi-homogeneous free divisors (e.g. free
hyperplane arrangements or discriminant of left-right stable
mappings in Mather's ``nice dimensions") are example of Koszul free
divisors \cite{calde_nar_LQHKF}.

\begin{definicion} 1) Let
$E\subset F$ be a pair of free $A$-modules of finite rank. We say
that $(E,F)$ is a {\em Koszul pair} (over $A$) if some (or any)
basis of $E$ forms a
regular sequence in the symmetric algebra $\Sim_A(F)$.\\
2) Similarly, we say that a pair $(\E,\F)$ of locally free $\scr A
$-modules of finite rank, with  $\E \subset \F$,  is a {\em Koszul
pair} (over $\scr A$) if $(\E_p,\F_p)$ is a Koszul pair over $({\scr
K}_p,{\scr L}_p)$ for any point $p\in M$.
\end{definicion}

To say that $(\Der(\log D), \Der_{\C}(\OX))$ is a Koszul pair is
equivalent to saying that $D$ is a Koszul free divisor.

\begin{definicion} 1) Let $L \subset L'$ be a pair of Lie-Rinehart algebras over
$(k,A)$.  We say that $(L,L')$ is a {\em pre-Spencer pair}
 (over $(k,A)$) if the complex $\U(L') \Lotimes_{\U(L)} A$
is cohomologically concentrated in degree $0$.\\
2) Similarly, we say that a pair $(\LL,\LL')$ of Lie algebroids over
$({\scr K},{\scr A})$  is a {\em pre-Spencer pair} if
$(\LL_p,\LL'_p)$ is a pre-Spencer pair over $({\scr K}_p,{\scr
A}_p)$ for any $p\in M$, or equivalently, if the complex
$\U(\LL')\Lotimes_{\U(\LL)} {\scr
A}$ is cohomologically concentrated in degree $0$.\\
3) We say that $D$ is a pre-Spencer (free) divisor if $(\Der(\log
D), \Der_{\C}(\OX))$ is a pre-Spencer pair over $(\C_X,\OX)$.
\end{definicion}

From (\ref{eq:notice}) and proposition \ref{prop:sp-res} we know that
$$ \SP_{L,L'}=\U(L')\otimes_{\U(L)} \SP_L(A) =\U(L') \Lotimes_{\U(L)} A,$$
and the property for $(L,L')$ to be a pre-Spencer pair is equivalent
to the fact that the complex $\SP_{L,L'}$ is cohomologically
concentrated in degree $0$, and so it is a free resolution of
$\U(L')/\U(L')\cdot L$ trough the augmentation (\ref{eq:nat-aug}).
In particular, if $D$ is a Spencer divisor in the sense of
\cite{cas_ucha_stek}, then it is a pre-Spencer divisor.

\begin{proposicion} \label{prop:K->Sp-E}
Let $L \subset L'$ be a pair of Lie-Rinehart algebras over $(k,A)$.
If $(L,L')$ is a Koszul pair over $A$ and $E$ is a left
$\U(L)$-module flat (resp. free) over $A$, then the complex
$\SP_{L,L'}(E)$ is a $\U(L')$-resolution (resp. a free
$\U(L')$-resolution) of $\U(L')\otimes_{\U(L)} E$.
\end{proposicion}

\begin{prueba} The proof is similar to the proof of proposition
\ref{prop:sp-res}. We consider the filtration on the complex
$\SP_{L,L'}(E)$ given by
$$F^i \SP^{-k}_{L,L'}(E) = (F^{i-k}\U(L')) \otimes_A
\stackrel{k}{\bigwedge} L \otimes_A E, i\geq 0.$$ Its graded complex
is, by using the Poincar\'{e}-Birkhoff-Witt theorem, canonically
isomorphic to the tensor product by $-\otimes_A E$ of the complex
$\Sim_A (L') \otimes_A \bigwedge^{\bullet} L$, where the
differential is given by
$$ d^{-k}
(P\otimes(\lambda_1\wedge\cdots\wedge\lambda_k)) =\sum_{i=1}^k
(-1)^{i-1}
(P\lambda_i)\otimes(\lambda_1\wedge\cdots\wedge\widehat{\lambda_i}
\wedge\cdots\wedge\lambda_k)$$ for $P\in \Sim_A (L')$,
$\lambda_1,\dots,\lambda_k\in L$ and $k=1,\dots, \rk_A L$.

Since $(L,L')$ is a Koszul pair and $E$ is flat over $A$, the
complex
$$ \Gr_F \SP_{L,L'}(E) = \left(\Sim_A (L') \otimes_A \bigwedge^{\bullet} L \right)
\otimes_A E$$ is exact in degrees $\neq 0$, and so $\SP_{L,L'}(E)$
too, i.e. it is a resolution of its $0$-cohomology $ h^0\left(
\SP_{L,L'}(E) \right) = \U(L')\otimes_{\U(L)} E$.
\end{prueba}

\begin{corolario} \label{cor:K->Sp} 1) Let $L \subset L'$ be a pair of Lie-Rinehart algebras
over $(k,A)$. If $(L,L')$ is a Koszul pair over $A$, then $(L,L')$
is
a pre-Spencer pair.\\
2) In a similar way, any pair $(\LL,\LL')$ of Lie algebroids over
$({\scr K},{\scr A})$ which is a Koszul pair over $\scr A$  is a
pre-Spencer pair.
\end{corolario}
\begin{prueba}  The second part follows straightforward from the first
part, and the first part is a consequence of Proposition
\ref{prop:K->Sp-E} in the case $E=A$.
\end{prueba}

\begin{proposicion} \label{prop:k[s]} Let $L\subset L'$ be a pair of $A$-modules (resp. of Lie-Rinehart algebras
over $(k,A)$) which are $A$-free  of finite rank. If $(L,L')$ is a
Koszul pair over $A$ (resp. a pre-Spencer pair over $(k,A)$), then
$(L[s],L'[s])$ is a Koszul pair over $A[s]$ (resp. a pre-Spencer
pair over $(k[s],A[s])$).
\end{proposicion}

\begin{prueba} It comes
from the flatness of $k \to k[s]$ and the fact that
$\Sim_{A[s]}(L[s]) = k[s]\otimes_k \Sim_A(L)$, $\bigwedge_{A[s]}
L[s] = k[s]\otimes_k \left( \bigwedge_A L \right) $
 and $$\U_{(k[s],A[s])}( L'[s]) = k[s]\otimes_k \U_{(k,A)}(L').$$
\end{prueba}

\subsection{The logarithmic Bernstein construction for free divisors}

In the situation of section \ref{subsec:D-bernstein}, let us suppose
that $D$ is a free divisor and let us write $\VO=\DX(\log D)_p$.
Since $D$ is free, the Lie-Rinehart algebra $\Thetafs$ defined in
(\ref{eq:thetafs}) is also $\hol$-free of rank $n$.

Similar to the case of $\D[s]$, the filtered ring $(\VO[s],F)$ is
the enveloping algebra of the Lie-Rinehart algebra $\Der(\log
D)_p[s]$ over $(\hol[s],\C[s])$, and the filtered ring
$(\VO[s],F_T)$ is the enveloping algebra of the Lie-Rinehart
algebra $F^1\VO= \hol\oplus \Der(\log D)_p$ over $(\hol,\C)$.

The free module of rank one over the ring $\hol[s]$ generated by the
symbol $f^s$, $\hol[s]f^s$, has a natural left module structure over
the ring $\VO[s]$: the action of a logarithmic derivation
$\delta\in\Der(\log D)_p$ is given by $\delta (f^s) = \frac{\delta
(f)}{f} s  f^s$.
\medskip

Let $\{\delta_1,\dots,\delta_n\}$ be a basis of $\Der(\log D)_p$ and
let us write $\delta_i(f)=\alpha_i f$ and $\eta_i=\sigma(\delta_i)\in
\Gr_F \VO = \hol[\eta_1,\dots,\eta_n]$ for $i=1,\dots,n$. The
$\zeta_i= \delta_i - \alpha_i s$, $i=1,\dots,n$, form a basis of
$\Thetafs$. Their symbols with respect to the total order filtration
are
$$\sigma_{T}(\zeta_i) = \eta_i - \alpha_i s \in
\Gr_{F_T}\VO[s]=\Sim_{\hol}(F^1\VO)= \hol[s,\eta_1,\dots,\eta_n],$$
and so $(\Thetafs,F^1\VO)$ is a Koszul pair. From corollary
\ref{cor:K->Sp} we deduce that the complex $\SP_{\Thetafs,F^1\VO}$
is cohomologically concentrated in degree $0$. On the other hand, a
division argument and the existence of the unique expressions
(\ref{eq:estruc-V0}) show that $\ann_{\VO[s]} f^s = \VO[s]
\Thetafs$. Finally, by using the augmentation (\ref{eq:nat-aug}) we
obtain a proof of the following proposition.

\begin{proposicion} \label{prop:VO-resol-fs}
Under the above conditions,
the complex $\SP_{\Thetafs,F^1\VO}$ is a free resolution of the
$\VO[s]$-module $\hol[s]f^s$.
\end{proposicion}

\begin{proposicion} \label{prop:thetafs-KF} If $D$ is Koszul free at $p$ then
$(\Thetafs,F^1\D)$ is a Koszul pair over $\hol$.
\end{proposicion}

\begin{prueba}
Let us take a system of local coordinates $x_1,\dots,x_n\in \hol$
and consider the symbols of the partial derivatives $\xi_i =
\sigma\left( \frac{\partial}{\partial x_i} \right)$. With the
notations above, we know that $\eta_1,\dots,\eta_n$ form a regular
sequence in $\Gr_F \D = \Sim_{\hol}(\Der_{\C}(\hol))=
\hol[\xi_1,\dots,\xi_n]$. Then, $s,\eta_1,\dots,\eta_n$ form another
regular sequence in $\hol[s,\xi_1,\dots,\xi_n] =
\Sim_{\hol}(F^1\D)$. Since the ideals $(s,\eta_1,\dots,\eta_n)$ and
$(s,\eta_1-\alpha_1 s,\dots,\eta_n-\alpha_n s)$ coincide, and we are
dealing with homogeneous elements in graded rings, regular and
quasi-regular sequence are the same (cf. \cite{mat_86}, \S~16) and
$s,\eta_1-\alpha_1 s,\dots,\eta_n-\alpha_n s$ is also a regular
sequence in $\hol[s,\uxi]$. In particular,
$\sigma_T(\zeta_1),\dots,\sigma_T(\zeta_n)$ is a regular sequence in
$\Sim_{\hol}(F^1\D)$ and the proposition is proved.
\end{prueba}

\begin{nota} Example \ref{ej:1} shows that, in the above proposition, the assumptions of being Koszul for $D$ is not
necessary to have the Koszul property for $(\Thetafs,F^1\D)$.
\end{nota}

The following theorem was announced in Remark 5.10, (a) in
\cite{calde_nar_compo}.

\begin{teorema} \label{teo:compo-5.10} Let us suppose that $D$ is of differential linear
type and Koszul free. Then the complex $\SP_{\Thetafs,F^1\D}$ is a
free resolution of the $\D[s]$-module $\D[s]\cdot f^s \subset
\hol[f^{-1},s]f^s$ and the canonical morphism
$$ \D[s] \Lotimes_{\VO[s]} \hol[s]f^s \xrightarrow{} \D[s]\cdot f^s$$
is an isomorphism.
\end{teorema}

\begin{prueba}  From Proposition \ref{prop:VO-resol-fs} and (\ref{eq:notice})
we have
$$ \D[s] \Lotimes_{\VO[s]} \hol[s]f^s = \D[s] \otimes_{\VO[s]}
\SP_{\Thetafs,F^1\VO}= \SP_{\Thetafs,F^1\D}. $$ On the other hand,
by Proposition \ref{prop:thetafs-KF}, the complex
$\SP_{\Thetafs,F^1\D}$ is exact in degrees $\neq 0$. To conclude, we
use that $D$ is of differential linear type:
$$ h^0\left(\SP_{\Thetafs,F^1\D}\right) = \D[s]/\D[s]\cdot \Thetafs
= \D[s]/\ann_{\D[s]} f^s = \D[s]\cdot f^s.$$
\end{prueba}

\begin{nota}  \label{nota:coment-LCT}
The hypotheses in the above theorem are satisfied in the case of
locally quasi-homogeneous free divisors. In fact, following the
lines in the proof of Theorem \ref{teo:main} and Theorem 4.1 in
\cite{calde_nar_fourier}, it is possible to deduce a weak (local)
version of the logarithmic comparison theorem (LCT). Namely, under
the hypothesis of theorem \ref{teo:compo-5.10}, there is a $k\gg 0$
such that the canonical morphism
$$ \D \Lotimes_{\VO} \left(\hol\cdot f^{-k}\right) \to
\hol[f^{-1}]$$ is an isomorphism in the derived category of
$\D$-modules (see \cite{calde_nar_compo}, Remark 5.10,~(b)). To go
further and deduce a proof of the full version of the LCT for a
locally quasi-homogeneous free divisors, one should prove before
that that the Bernstein-Sato polynomial of its (reduced) equation
does not have any integer root less than $-1$. Unfortunately, we do
not know any direct proof of this fact. Nevertheless, see
\cite{calde_nar_fourier}, Th.~4.4 and Corollary \ref{cor:LCT} for
other proofs of the LCT based on $D$-module theory.
\end{nota}

The following remark has been pointed out by Torrelli (part [a] was
also known by the authors).
\medskip

\begin{nota} \label{nota:torre}
Let $p$ be a point in $D$, $x_1,\dots,x_n\in{\scr O}_{X,p}$ a system
of local coordinates centered at $p$ and $f\in\hol=\hol_{X,p}$ a
reduced local equation of $D$.
\medskip

\noindent [a] We know that $f$ belongs to the integral closure of
the gradient ideal $I=(f'_{x_1},\dots,f'_{x_n})$ (cf.
\cite{MR0374482}, \S 0.5, 1)), i.e there is an integer $d>0$ and
elements $a_i \in I^{d-i}$ such that $f^d + a_{d-1} f^{d-1} + \cdots
+ a_0 =0$. In particular, there is a homogeneous polynomial $F\in
\hol[s,\xi_1,\dots,\xi_n]$ of degree $d>0$ such that
$F(f,f'_{x_1},\dots,f'_{x_n})=0$ and $F(s,0,\dots,0)=s^d$. Let $$
\delta_i = \sum_{j=1}^n a_{ij} \frac{\partial }{\partial x_j},\quad
1\leq i\leq m$$ a system of generators of $\Der(\log D)_p$ and let
us write $\delta_i(f)=\alpha_i f$. In other words,
$(-\alpha_i,a_{i1},\dots,a_{in})$, $1\leq i\leq m$, is a system of
generators of the syzygies of $f,f'_{x_1},\dots,f'_{x_n}$. If $D$ is
of linear jacobian type at $p$, then the polynomial $F$ must be a
linear combination of the polynomials
$$ -\alpha_i s + a_{i1} \xi_1 + \cdots +a_{in}\xi_n, \quad 1\leq i\leq m,$$
and making $\xi_1=\cdots=\xi_n=0$ we deduce that some of the
$\alpha_i$ must be a unit, i.e. $f \in (f'_{x_1},\dots,f'_{x_n})$.
That shows that if $D=\{f=0\}$ is of linear jacobian type then $f$
is {\em Euler homogeneous}, i.e. there is a germ of vector field
$\chi$ such that $\chi(f)=f$.
\medskip

\noindent [b] Assume that the annihilator of $f^{-1}$ over $\D$ is
generated by operators of order one. Then, from Proposition 1.3 of
\cite{torre-45-bis} we know that $-1$ is the smallest integer root
of the Bernstein polynomial $b_f(s)$ of $f$. Reciprocally, let us
assume that $D$ is of linear jacobian type at $p$ and that $-1$ is
the smallest integer root of $b_f(s)$. Then, the annihilator of
$f^{-1}$ over $\D$ is generated by operators of order one since it
is obtained from the annihilator of $f^s$ over $\D[s]$ by making
$s=-1$.
\medskip

\noindent [c] Let us suppose now that $D$ is a free divisor which is
Koszul and of differential linear type at $p$, and let $m_0$ be the
smallest integer root of the Bernstein polynomial of $f$. For $l\geq
-m_0$ the annihilator of $f^{-l}$ over $\D$ is generated by
operators of order one, and the Koszul hypothesis allows us to apply
Proposition 4.1 of \cite{torre-45-bis} in order to obtain that $f$
is Euler homogeneous.
\end{nota}

The following result is proved in \cite{simis-alg-free}, Cor.~3.12
in the polynomial case and has been independently pointed out to us
 by T. Torrelli.

\begin{proposicion} \label{prop:simis-torr}  If $D$ is of
linear jacobian type and free at $p$, then it is Koszul free at $p$.
\end{proposicion}

\begin{prueba} Let $x_1,\dots,x_n\in{\scr O}_{X,p}$ be a system
of local coordinates centered at $p$ and $f\in\hol=\hol_{X,p}$ a
reduced local equation of $D$. From remark \ref{nota:torre}, a) we
know that $f$ is Euler homogeneous, i.e. $f\in
(f'_{x_1},\dots,f'_{x_n})$ and so $\jac(D)_p
=(f'_{x_1},\dots,f'_{x_n})$.

Let $\{\delta_i = \sum_{j=1}^n a_{ij} \frac{\partial}{\partial
x_j}\}_{1\leq i \leq n}$ be a basis of $\Der(\log D)_p$ and let us
write $\delta_i(f)=\alpha_i f$. Since $f$ is Euler homogeneous, we
can take $\alpha_1=\cdots=\alpha_{n-1}=0$ and $\alpha_n=1$. In other
words, $\{(a_{i1},\dots,a_{in}\}_{1\leq i \leq n-1})$ is a basis of
the syzygies of $f'_{x_1},\dots,f'_{x_n}$.

Let $\theta : \hol[\xi_1,\dots,\xi_n] \xrightarrow{}
\Rees(\jac(D)_p)=\hol[f'_{x_1}t,\dots,f'_{x_n}t]$ be the surjective
map of $\hol$-algebras defined by $\theta(\xi_i)=f'_{x_i}t$. Since
$\jac(D)_p$ is an ideal of linear type, the kernel of $\theta$ is
generated by the $\sigma(\delta_i)= \sum_{j=1}^n a_{ij}\xi_j$,
$1\leq i\leq n-1$. So
$$\dim \left(
\frac{\hol[\xi_1,\dots,\xi_n]}{(\sigma(\delta_1),\dots,\sigma(\delta_{n-1}))}\right)
=\dim \Rees(\jac(D)_p) = n+1$$ and
$\sigma(\delta_1),\dots,\sigma(\delta_{n-1})$ is a regular sequence.

On the other hand if $F\sigma(\delta_n) \in
(\sigma(\delta_1),\dots,\sigma(\delta_{n-1}))$, so $0=
\theta(F)\theta(\sigma(\delta_n))= \theta(F) ft$ and we deduce that
 $F\in \ker \theta$ and
$\sigma(\delta_1),\dots,\sigma(\delta_n)$ is a regular sequence.
\end{prueba}

\begin{ejemplo} \label{ej:1}
Let us suppose that $D\subset X$ is a non-necessarily free divisor
and let $f=0$ be a reduced local equation of $D$ at a point $p\in
D$. Let $\{\delta_1,\dots,\delta_m\}$ be a system of generators of
$\Der(\log D)_p$ and let us write $\delta_i(f)=\alpha_i f$.

Let us call $\ann_{\D[s]}^{(1)}(f^s)$ the ideal of $\D[s]$ generated
by $\Thetafs$ (\ref{eq:thetafs}):
$$ \ann_{\D[s]}^{(1)}(f^s) = \D[s]\cdot\left(\delta_1 - \alpha_1
s,\dots,\delta_m -\alpha_m s\right)\subset \ann_{\D[s]}(f^s).$$

The Bernstein functional equation for $f$ \cite{bern_72,kas_76}
$$ b(s) f^s = P(s) f^{s+1}$$ means that the operator $b(s) - P(s) f$
belongs to the annihilator of $f^s$ over $\D[s]$. Then, an explicit
knowledge of the ideal $\ann_{\D[s]}(f^s)$ allows us to find $b(s)$
by computing the ideal $\C[s]\cap \left( \D[s]\cdot f +
\ann_{\D[s]}(f^s) \right)$. However, the ideal $\ann_{\D[s]}(f^s)$
is in general difficult to compute.

When $D$ is a divisor of differential linear type,
$\ann_{\D[s]}(f^s)= \ann_{\D[s]}^{(1)}(f^s)$ and the computation of
$b(s)$ is in principle easier. But there are examples of free
divisors which are not of differential linear type for which the
Bernstein polynomial $b(s)$ belongs to
$$\C[s]\cap \left( \D[s]\cdot f +
\ann_{\D[s]}^{(1)}(f^s) \right).$$ For instance, for $X=\C^3$ and
$f=x_1 x_2 (x_1+x_2) (x_1+x_2 x_3)$ (see Example 6.2 in
\cite{calde_nar_compo}) or in the examples in page 445 of
\cite{cas_ucha_exper}. In all this examples the divisor is not
Koszul, satisfies the logarithmic comparison theorem and
$(\Thetafs,F^1\D)$ is a Koszul pair over $\hol$ (see
Prop.~\ref{prop:thetafs-KF}).
\end{ejemplo}

\section{Integrable logarithmic connections with res\-pect to a free divisor}

In this section we assume that $D\subset X$ is a free divisor.

A {\em logarithmic connection} with respect to $D$ is a locally
$\OX$-module $\E$ (the case where $\E$ is only supposed to be
coherent certainly deserves
to be studied, but it will not be treated in this paper) endowed with:\\
-) a $\C$-linear map (connection) $ \nabla': \E\xrightarrow{}
\E\otimes_{\OX}\Omega^1_X(\log D)$ satisfying $\nabla'(ae) = a
\nabla'(e) + e\otimes da$, for any section $a$
of $\OX$ and any section $e$ of $\E$,\\
or equivalently, with\\
-) a left $\OX$-linear map $\nabla:\derlogD \xrightarrow{}
\End_{\C_X}(\E)$ satisfying the Leibniz rule $\nabla(\delta)(ae) =
a\nabla(\delta)(e) + \delta(a)e$, for any logarithmic vector field
$\delta$, any section $a$ of $\OX$ and any section $e$ of $\E$.

The integrability of $\nabla'$ is equivalent to the fact that
$\nabla$ preserves Lie brackets. Then, we know from
\ref{nume:carac-Umod} and section \ref{subsec:log-der} that giving
an integrable logarithmic connection on a locally free $\OX$-module
$\E$ is equivalent to extending its original $\OX$-module structure
to a left $\DX(\log D)$-module structure, and so integrable
logarithmic connections are the same as left $\DX(\log D)$-modules
which are locally free of finite rank over $\OX$.

Let us denote by $\OX(\star D)$ the sheaf of meromorphic functions
with poles along $D$. It is a holonomic left $\DX$-module
\cite{kas_76}.

The first examples of integrable logarithmic connections (ILC for
short) are the invertible $\OX$-modules $\OX(mD)\subset \OX(\star
D)$, $m\in\ZZ$, formed by the meromorphic functions $h$ such that
$\divi (h) + mD \geq 0$.

If $f=0$ is a reduced local equation of $D$ at $p\in D$ and
$\delta_1,\dots,\delta_n$ is a local basis of $\Der(\log D)_p$ with
$\delta_i(f)=\alpha_i f$, so $f^{-m}$ is a local basis of
$\OO_{X,p}(mD)$ over $\OO_{X,p}$ and we have the following local
presentation over $\D_{X,p}(\log D)$ (use (\ref{eq:estruc-V0}))
\begin{equation*}
\OO_{X,p}(mD) \simeq \D_{X,p}(\log D)/\D_{X,p}(\log
D)(\delta_1+m\alpha_1,\dots, \delta_n+m\alpha_n).\end{equation*}

For any ILC $\E$ and any integer $m$, the locally free $\OX$-modules
$\E(mD):= \E \otimes_{\OX} \OX(mD)$ and $\E^* := \Hom_{\OX}(\E,\OX)$
are endowed with a natural structure of left $\DX(\log D)$-module
(cf. \cite{calde_nar_fourier}, \S 2), and they are again ILC, and
the usual isomorphisms
$$ \E(mD)(m'D) \simeq \E((m+m')D),\quad \E(mD)^* \simeq \E^*(-mD)$$are
$\DX(\log D)$-linear.

\subsection{The logarithmic comparison problem}
\label{subsec:LCP}

If $D$ is Koszul free and $\E$ is an integrable logarithmic
connection, then the complex $\DX\Lotimes_{\DX(\log D)} \E$ is
concentrated in degree $0$ and its $0$-cohomology
$\DX\otimes_{\DX(\log D)} \E$ is a holonomic $\DX$-module (see
\cite{calde_nar_fourier}, Prop.~1.2.3).

Let us denote by $\DX(\star D)$ the sheaf of meromorphic linear
differential operators with poles along $D$. One has obvious left and
right $\OX(\star D)$-linear isomorphisms
$$ \OX(\star D)\otimes_{\OX} \DX \stackrel{\text{left}}{\simeq}
\DX(\star D) \stackrel{\text{right}}{\simeq} \DX \otimes_{\OX}
\OX(\star D).$$The induced maps
$$ \OX(\star D)\otimes_{\OX} \DX(\log D) \xrightarrow{} \DX(\star D)
\xleftarrow{} \DX(\log D) \otimes_{\OX} \OX(\star D)$$are also
isomorphisms and so ``meromorphic logarithmic linear differential
operators" and ``meromorphic linear differential operators" are the
same:
$$ \DX(\log D)(\star D) = \DX(\star D).$$

If $\E$ is a left $\DX(\log D)$-module, then the localization
\begin{equation} \label{eq:locali}
\E(\star D) := \OX(\star D)\otimes_{\OX} \E = \DX(\star
D)\otimes_{\DX(\log D)} \E
\end{equation} is a left $\DX(\star
D)$-module, and by scalar restriction, a left $\DX$-module.
Moreover, if $\E$ is a ILC, then $\E(\star D)$ is a meromorphic
connection (locally free of finite rank over $\OX(\star D)$) and
then it is a holonomic $\DX$-module (cf. \cite{meb_nar_dmw}, Th.
4.1.3). Actually, $\E(\star D)$ has regular singularities on the
smooth part of $D$ (it has logarithmic poles! \cite{del_70}) and so
it is regular everywhere \cite{meb-cimpa-2}, Cor. 4.3-14, which
means that if $\LL$ is the local system of horizontal sections of
$\E$ on $U=X-D$, the canonical morphism $
\Omega_X^{\bullet}(\E(\star D)) \xrightarrow{} R j_* \LL$ is an
isomorphism in the derived category.

For any ILC $\E$, or even for any left $\DX(\log D)$-module (without
any finiteness property over $\OX$), one can define its logarithmic
de Rham complex $\Omega^{\bullet}_X(\log D)(\E)$ in the classical
way, which is a subcomplex of $\Omega^{\bullet}_X(\E(\star D))$. It
is clear that both complexes coincide on $U$.

For any ILC $\E$ and any integer $m$, $\E(mD)$ is a sub-$\DX(\log
D)$-module of the regular meromorphic connection (and holonomic
$\DX$-module) $\E(\star D)$, and so we have a canonical morphism in
the derived category of left $\DX$-modules
\begin{equation} \label{eq:rhoEm}
 \rho_{\E,m}: \DX \Lotimes_{\DX(\log D)} \E(mD)\to \E(\star D),\end{equation}
given by $\rho_{\E,m}(P\otimes e') = Pe'$.

Since $\E(m'D)(mD) = \E((m+m')D)$ and $\E(m'D)(\star D) = \E(\star
D)$, we can identify morphisms $\rho_{\E(m'D),m}$ and
$\rho_{\E,m+m'}$.

We have the following theorem:

\begin{teorema} \label{teo:crit-LCT}
Let $\E$ be a ILC (with respect to the free divisor $D$) and let
$\LL$ be the local system of its horizontal sections on $U=X-D$. The
following properties are equivalent:
\begin{enumerate}
\item[1)] The canonical morphism $\Omega_X^{\bullet}(\log D)(\E)
\to R j_* \LL$ is an isomorphism in the derived category of
complexes of sheaves of complex vector spaces.
\item[2)] The inclusion $\Omega_X^{\bullet}(\log D)(\E) \hookrightarrow
\Omega_X^{\bullet}(\E(\star D))$ is a quasi-isomorphism.
\item[3)] The morphism $\rho_{\E,1}: \DX\Lotimes_{\DX(\log D)} \E(D)\to \E(\star D)$
is an isomorphism in the derived category of left $\DX$-modules.
\item[4)] The complex $\DX\Lotimes_{\DX(\log D)} \E(D)$ is
concentrated in degree $0$ and the $\DX$-module
$\DX\otimes_{\DX(\log D)} \E(D)$ is holonomic and isomorphic to its
localization along $D$.
\item[5)] The canonical morphism $j_!\LL^{\vee} \to
\Omega_X^{\bullet}(\log D)(\E^*(-D))$ is an isomorphism in the
derived category of complexes of sheaves of complex vector spaces.
\end{enumerate}
\end{teorema}

\begin{prueba} The equivalence of the first three
properties has been proved in \cite{calde_nar_fourier}, Th. 4.1.

The equivalence between 3) and 4) comes from the fact that the
localization along $D$ of $\DX\Lotimes_{\DX(\log D)} \E(D)$ is
canonically isomorphic to $\E(\star D)$:
\begin{eqnarray*}
& \left[ \DX\Lotimes_{\DX(\log D)} \E(D) \right](\star D) \simeq
\DX(\star D) \otimes_{\DX}  \left[ \DX\Lotimes_{\DX(\log D)} \E(D)
\right]\simeq&\\
& \DX(\star D) \otimes_{\DX(\log D)} \E(D)\simeq \OX(\star D)
\otimes_{\OX} \E(D)\simeq \E(D)(\star D) \simeq \E(\star D).&
\end{eqnarray*}
The equivalence between 5) and 1) is a consequence of the duality
result in \cite{calde_nar_fourier}, Cor.~3.1.8,
$$ \Omega_X^{\bullet}(\log D)(\E)^{\vee} \simeq \Omega_X^{\bullet}(\log
D)(\E^*(-D)),$$ \cite{calde_nar_fourier}, Cor.~3.1.6 and the fact
that $ \left( R j_* \LL\right)^{\vee} = j_!\LL^{\vee}$.
\end{prueba}

\begin{nota} In the above theorem, the complex $\DX\Lotimes_{\DX(\log D)}
\E(D)$ need not to be holonomic even if its localization along $D$
is holonomic. For instance, in Example 5.1 of
\cite{calde_nar_fourier} for $\E=\OX(-D)$ the $\DX$-module
$\DX\otimes_{\DX(\log D)} \OX$ is not holonomic. This fact has been
also pointed out by Castro-Ucha and Torrelli. In particular, Problem
5.4 in {\em loc.~cit.} has a negative answer.
\end{nota}

For $D$ a locally quasi-homogeneous free divisor and $\E=\OX$, the
equivalent properties in Theorem \ref{teo:crit-LCT} hold: this is
the so called ``logarithmic comparison theorem"
\cite{cas_mond_nar_96}, \cite{calde_nar_fourier}, Th. 4.4. Here we
give a new proof using property 5) in Theorem \ref{teo:crit-LCT}.

\begin{proposicion} \label{prop:mond} Let us suppose that $D$ is a locally quasi-homogeneous (not necessarily free) divisor.
Then the canonical morphism $$j_! \C_U \to \Omega_X^{\bullet}(\log
D)(\OX(-D))$$ is a quasi-isomorphism.
\end{proposicion}

\begin{prueba} By Poincar\'e's lemma, the result is clear on $U$. To
conclude, we apply \cite{mond-plms-2000}, Lemma 3.3, (6) and we
deduce that the complex $\Omega_X^{\bullet}(\log D)(\OX(-D))$ is
acyclic at any point $p\in D$.
\end{prueba}

\begin{corolario} \cite{cas_mond_nar_96} \label{cor:LCT} Let $D$ be a locally quasi-homogeneous
free divisor. Then the {\em logarithmic comparison theorem} holds:
$$\Omega_X^{\bullet}(\log D) \xrightarrow{\sim} R j_* \C_U.$$
\end{corolario}

\begin{prueba} The result is a straightforward
consequence of Theorem \ref{teo:crit-LCT} and Proposition
\ref{prop:mond}.
\end{prueba}

\section{Main results}

Throughout this section, we suppose that $D\subset X$ is a free
divisor and $\E$ is an ILC with respect to $D$.

\subsection{The Bernstein-Kashiwara construction for integrable logarithmic
connections} \label{subsec:BK}

Let $p$ be a point in $D$ and $f\in\hol=\hol_{X,p}$ a reduced local
equation of $D$. Let us write $\D= \D_{X,p}$, $\VO=\DX(\log D)_p$,
$\Der(\log f)=\Der(\log D)_p$ and $E=\E_p$.

We know that $\E(\star D)$ (resp. $E[f^{-1}]= \E(\star D)_p$) is a
left $\DX$-module (resp. a left $\D$-module) (see
(\ref{eq:locali})).

The module $E[f^{-1},s]f^s = E\otimes_{\hol} \hol[f^{-1},s]f^s =
E[f^{-1}]\otimes_{\hol} \hol[f^{-1},s]f^s$ has a natural module
structure over the ring $\D[s]$: the action of a derivation
$\delta\in\Der_{\C}(\hol)$ is given by $\delta(e f^s) = \delta(e)
f^s + s\delta(f) f^{-1}  e f^s$. We have $\VO[s]\cdot E f^s =
E[s]f^s$, and so $E[s]f^s$ is a sub-$\VO[s]$-module of
$E[f^{-1},s]f^s$.

From Proposition \ref{prop:sp-res} we know that the complex
$\SP_{\Der(\log f)[s]}(E[s]f^s)$ is a free $\VO[s]$-resolution of
$E[s]f^s$ (here we consider $\VO[s]$ as the enveloping algebra of
the Lie-Rinehart algebra $\Der(\log f)[s]$ over $(\C[s],\hol[s])$).
On the other hand, we have a canonical $\D[s]$-linear map
$$ P\otimes (ef^s) \in \D[s]\otimes_{\VO[s]} E[s]f^s \mapsto P(ef^s)
\in \D[s]\cdot E[s]f^s \subset E[f^{-1},s]f^s$$ inducing a
surjective augmentation \begin{equation} \label{eq:augEfs}
\rho_{E,s}:\D[s]\otimes_{\VO[s]} \SP^0_{\Der(\log f)[s]}(E[s]f^s)
\xrightarrow{} \D[s]\cdot \left(E[s]f^s\right) = \D[s]\cdot \left(Ef^s\right).\end{equation}

The following theorem is strongly related to Theorem
\ref{teo:compo-5.10} in the case of the trivial ILC $\E=\OX$.

\begin{teorema} \label{teo:main-0} Let us suppose that $D$ is of
linear jacobian type at $p\in D$. Then the complex
$$\D[s]\otimes_{\VO[s]} \SP_{\Der(\log f)[s]}(E[s]f^s)$$is exact and
becomes a free $\D[s]$-resolution of $\D[s]\cdot \left(Ef^s\right)$ through the
map $\rho_{E,s}$ in (\ref{eq:augEfs}).
\end{teorema}

\begin{prueba}
From (\ref{eq:notice}) and Propositions \ref{prop:K->Sp-E},
\ref{prop:simis-torr} and \ref{prop:k[s]}, we deduce that  the
complex $\D[s]\otimes_{\VO[s]} \SP_{\Der(\log f)[s]}(E[s]f^s)$ is
exact in degrees $\neq 0$. To conclude, we need to prove that the
sequence
$$\D[s]\otimes_{\hol[s]} \Der(\log f)[s] \otimes_{\hol[s]} E[s]f^s \xrightarrow{\varepsilon_s^{-1}}
\D[s] \otimes_{\hol[s]} E[s]f^s \xrightarrow{\rho_{E,s}} \D[s]\cdot
\left(Ef^s\right)$$is exact, where $\varepsilon_s^{-1}(P\otimes
\delta \otimes (e f^s)) = (P\delta)\otimes (e f^s) - P\otimes
\delta(e f^s)$ (see \ref{nume:CECRSP}) and $\rho_{E,s}(P\otimes (e
f^s)) = P(e f^s)$. The inclusion $\im \varepsilon_s^{-1} \subset
\ker \rho_{E,s}$ is clear.

Let $\{e_1,\dots,e_r\}$ be an $\hol$-basis of $E$. Any $Q\in \D[s]
\otimes_{\hol[s]} E[s]f^s $ can be uniquely written as $Q =
\sum_{i=1}^r Q_i \otimes e_i f^s$ with $Q_i\in\D[s]$. We define the
total order of $Q$, $\Deg_T(Q)$, as the maximum of the orders of the
$Q_i$ with respect to the total order filtration in $\D[s]$ and
$$F_T^k\left( \D[s] \otimes_{\hol[s]} E[s]f^s \right) = \{Q\ |\ \Deg_T(Q) \leq k\}.$$

Let $Q = \sum_{i=1}^r Q_i \otimes e_i f^s \in \ker \rho_{E,s}$. To
prove that $Q$ belongs to the image of $\varepsilon_s^{-1}$, we
proceed by induction on $\Deg_T(Q)$. If $\Deg_T(Q)=0$, then the $Q_i$
belong to $\hol$ and the result is clear.  Let us suppose now that
$$F_T^{k-1}\left( \D[s] \otimes_{\hol[s]} E[s]f^s \right)\cap \ker \rho_{E,s}\subset
\im \varepsilon_s^{-1}$$ and $\Deg_T(Q)= k$. We have
$$ 0 = \rho_{E,s} (Q) = \left(\sum_{i\in I} C_{Q_i,k} f^{-k}e_i f^s   \right) s^k +\
\text{terms of lower degree in $s$},$$ where $I=\{i\ |\
\Deg_T(Q_i)=k\}$ and $\varphi(\sigma_T(Q_i))=C_{Q_i,k}t^k$ (see
Lemma \ref{lema:coro}). Consequently, $\sigma_T(Q_i) \in \ker
\varphi$ for any $i\in I$ and, from Proposition \ref{prop:clt->dlt}
and Remark \ref{nota:expli-clt}, there are $P_{ij} \in F^{k-1}_T
\D[s]$ and $\gamma_{ij}=\delta_{ij} - \alpha_{ij} s\in \Thetafs$
such that
$$ \sigma_T(Q_i) = \sum_j \sigma_T(P_{ij}) \sigma_T(\gamma_{ij}),\quad \forall i\in I.$$
Let us consider
$$ Q' = \sum_{i\in I} \left( \sum_j P_{ij}\otimes \gamma_{ij} \right) \otimes (e_i f^s) \in
\D[s]\otimes_{\hol[s]} \Der(\log f)[s] \otimes_{\hol[s]} E[s]f^s.$$
Since
\begin{eqnarray*} &\displaystyle \varepsilon_s^{-1} (Q')
= \sum_{i\in I}\left( \sum_j (P_{ij} \gamma_{ij}) \right)\otimes (e_i f^s) -
\sum_{i\in I}\left( \sum_j P_{ij}  \otimes \gamma_{ij}(e_i f^s)\right) =&\\
 &\displaystyle  \sum_{i\in I}\left( \sum_j (P_{ij} \gamma_{ij}) \right)\otimes (e_i f^s) -
\sum_{i\in I}\left( \sum_j P_{ij}  \otimes (\delta_{ij}\cdot e_i) f^s\right), &
\end{eqnarray*}
we have that $Q - \varepsilon_s^{-1} (Q') \in F_T^{k-1}\left( \D[s]
\otimes_{\hol[s]} E[s]f^s \right)\cap \ker \rho_{E,s}$, and by the
induction hypothesis we obtain that $Q$ belongs to $\im
\varepsilon_s^{-1}$.
\end{prueba}

\begin{corolario} \label{cor:main-0}
Under the hypothesis of Theorem \ref{teo:main-0}, the canonical
morphism
$$ \D[s]\Lotimes_{\VO[s]} E[s]f^s \to \D[s]\cdot (E f^s)$$
is an isomorphism in the derived category of left $\D[s]$-modules.
\end{corolario}

\begin{prueba} It is a consequence of Proposition \ref{prop:sp-res}
and Theorem \ref{teo:main-0}.
\end{prueba}

\begin{corolario} Let us suppose that $D$ is a locally quasi-homogeneous free divisor. Then the complex
$$\D[s]\otimes_{\VO[s]} \SP_{\Der(\log f)[s]}(E[s]f^s)$$is a free $\D[s]$-resolution of $\D[s]\cdot \left(Ef^s\right)$.
\end{corolario}

\begin{prueba} It is a consequence of Theorems \ref{teo:main-0} and \ref{teo:lqh->clt}, and
\cite{calde_nar_LQHKF}.
\end{prueba}

\subsection{The logarithmic comparison theorem}

Let us keep the notations of section \ref{subsec:BK}. Let us also
write $E[f^{-1}]= \E(\star D)_p$, $f^{-m}E = \E(mD)_p$ and
$$\rho_{E,m}: \D \otimes_{\hol} (f^{-m}E)
\xrightarrow{} E[f^{-1}]$$  the induced map by $\rho_{\E,m}$ in
(\ref{eq:rhoEm}), for any integer $m$.

\begin{lema}
There exists a non zero polynomial $b(s)\in\C[s]$ such that
$$ b(s) Ef^s \subset \D[s]\cdot \left(Ef^{s+1}\right).$$
\end{lema}

\begin{prueba} Let $\{e_1,\dots,e_r\}$ be
an $\hol$-basis of $E$ and let $b_i(s)$ be the Bernstein-Sato
polynomial of $e_i$ considered as an element of the holonomic $\D$-module $E[f^{-1}]$. We take
$b(s)=\lcm(b_1(s),\ldots,b_r(s))$.
\end{prueba}

\begin{nota}\label{nota:itera-BS}
a) The set of polynomials $b(s)$ in the above lemma is an ideal of
$\C[s]$, whose monic generator will be denoted by $b_E(s)$ (or
$b_{\E,p}(s)$) and will be called {\em Bernstein-Sato polynomial of
$\E$ at $p$}. \\
b) For any integer $k$ it is clear that the polynomials $b_E(s-k)$
satisfies $b_E(s-k) Ef^{s-k} \subset \D[s] \left(Ef^{s-k+1}\right)$.
In other words, $b_{\E,p}(s-k)=b_{\E(kD),p}(s)$.
 So, the polynomial $b_l(s)=\prod_{k=1}^l b_E(s-k)$
satisfies $b_l(s) Ef^{s-l} \subset \D[s]\left( Ef^{s}\right)$ and
$E[f^{-1}] = \D\cdot (f^{-m}E)$, i.e. $\rho_{E,m}$ is surjective, if
$b_E(s)$ has no integer roots less than $-m$.
\end{nota}

The following proposition can be useful in order to compute the polynomial $b_E(s)$.

\begin{proposicion} Let us suppose that $E$ is a cyclic $\VO$-module
generated by an element $e\in E$. Then, the polynomial $b_E(s)$
coincides with the Bernstein-Sato polynomial $b_e(s)$ of $e$ with respect to
$f$, where $e$ is considered as an element of the holonomic
$\D$-module $E[f^{-1}]$.
\end{proposicion}

\begin{prueba}
For any $\delta\in\Der(\log f)$, we have
$$ \delta (e f^{s+1}) = (\delta e) f^{s+1} + (s+1) \frac{\delta(f)}{f} e f^{s+1} $$
and so, since $\VO=\hol[\Der(\log f)]$, we deduce that $(\VO\cdot e)
f^{s+1} \subset \VO[s]\cdot (e f^{s+1})$ and
\begin{eqnarray*} &
 b_E(s) e f^s \in \D[s]\cdot \left(E f^{s+1}\right) = \D[s]\cdot\left( (\VO\cdot e) f^{s+1}\right) \subset&\\
 & \subset \D[s]\cdot\left( \VO[s]\cdot (e f^{s+1}) \right) = \D[s]\cdot \left(e f^{s+1}\right).\end{eqnarray*}
In particular $b_e(s)\ |\ b_E(s)$. On the other hand,
\begin{eqnarray*} &
b_e(s) E f^s = b_e(s) \left( (\VO\cdot e) f^s \right) \subset b_e(s) \left( \VO[s]\cdot (e f^s)\right) =&\\
&
\VO[s]\cdot \left( b_e(s) e f^s\right) \subset \VO[s]\cdot \left( \D[s]\cdot (e f^{s+1})\right) \subset \D[s]\cdot
\left( E f^{s+1}\right)&
\end{eqnarray*}
and so $b_E(s)\ |\ b_e(s)$.
\end{prueba}

\numero {\em Specialization at integers}
\medskip

For each $m\in \ZZ$ and $r\geq 0$, let us denote by $\Phi_m:\D[s] \to
\D$, $\Phi_{E,m}:E[f^{-1},s]f^s\to E[f^{-1}]$ and
$$\textstyle \Phi^r_{E,m}: \D[s]\otimes_{\hol[s]} \left(\bigwedge^r
\Der(\log f)[s]\right)\otimes_{\hol[s]} E[s]f^s \xrightarrow{} \D
\otimes_{\hol} \left(\bigwedge^r \Der(\log f)\right)\otimes_{\hol}
(f^{-m} E)$$ the specialization maps making $s=-m$. It is clear that
for any integer $m$ the following properties hold:
\begin{enumerate}
\item[-)]  $\Phi_{E,m}(E[s]f^s)\subset f^{-m} E$,
\item[-)] $\Phi_{E,m}$ and the $\Phi^r_{E,m}$ are $\Phi_m$-linear,
\item[-)] $\Phi_{E,m}\pcirc \rho_{E,s}=\rho_{E,m}\pcirc\Phi_{E,m}^0$,
\item[-)] the $\Phi^r_{E,m}$, $r\geq 0$, commute with the
differentials and define a morphism between
Cartan-Eilenberg-Chevalley-Rinehart-Spencer complexes (see
\ref{nume:CECRSP}).
\end{enumerate}

\begin{proposicion} \label{prop:kerrhos-kerrhok} Under the above conditions,
we have $\Phi_{E,k}^0(\ker \rho_{E,s}) = \ker \rho_{E,k}$ for all $k
\geq -m_0$, where $m_0$ is the smallest integer root of $b_E(s)$.
\end{proposicion}

\begin{prueba}
 Since $\Phi_{E,k}\pcirc \rho_{E,s}=\rho_{E,k}\pcirc\Phi_{E,k}^0$ for any integer $k$,
we deduce that $\Phi_{E,k}^0(\ker \rho_{E,s}) \subset \ker
\rho_{E,k}.$
\smallskip

\noindent Let $\{e_1,\dots,e_r\}$ be an $\hol$-basis of $E$ and
$P=\sum_i P_i\otimes e_i f^{-k}\in \ker \rho_{E,k}$. Let us consider
$P'= \sum_i P_i\otimes e_i f^s\in \D[s]\otimes_{\hol[s]} E[s]f^s$.
Since $\Phi_{E,k} (\rho_{E,s}(P')) = \rho_{E,k} (\Phi_{E,k}^0 (P'))
=\rho_{E,k}(P)=0$ we deduce that $\rho_{E,s}(P')$ is divisible by
$s+k$ in $E[f^{-1},s]f^s$, i.e. there is a $B\in E[s]f^s$ and an
$l>0$ such that $\rho_{E,s}(P') = (s+k) f^{-l} B$.

From Remark \ref{nota:itera-BS}, b), we know that $b_l(s) f^{-l} B
\in \D[s]\left(E f^s\right)$ and so there is a $Q\in
\D[s]\otimes_{\hol[s]} E[s]f^s$ such that $ b_l(s) f^{-l} B =
\rho_{E,s}(Q)$. The element $R=b_l(s) P' -(s+k)Q$ clearly belongs to
$\ker \rho_{E,s}$. If $k\geq -m_0$, $b_l(-k)\neq 0$ and
$$ P = \Phi_{E,k}^0\left(b_l(-k)^{-1}R\right) \in \Phi_{E,k}^0(\ker \rho_{E,s}).$$
\end{prueba}

\begin{teorema}  \label{teo:main} Let us suppose that $D$ is of
linear jacobian type and let $m_p$ be the smallest integer root of
$b_{\E,p}(s)$. Then, there is an open neighborhood $V$ of $p$ such
that the restriction to $V$ of the
 morphism $$ \rho_{\E,k}: \DX\Lotimes_{\DX(\log D)}
\E(kD) \xrightarrow{} \E(\star D)$$ is an isomorphism in the derived
category of $\D_U$-modules, for all $k\geq -m_p$.
\end{teorema}

\begin{prueba} By the coherence of the involved objects, we can work at the level of the
stalks at $p$.

 Since $D$ is Koszul (Prop. \ref{prop:simis-torr}), the complex
$$ \left[ \DX\Lotimes_{\DX(\log D)} \E(kD)\right]_p =
\D\Lotimes_{\VO} f^{-k}E$$is exact in degrees $\neq 0$ (see
Proposition \ref{prop:K->Sp-E}). To conclude, we need to prove that
the sequence
$$\D\otimes_{\hol} \Der(\log f) \otimes_{\hol} f^{-k}E \xrightarrow{\varepsilon_k^{-1}}
\D \otimes_{\hol} f^{-k}E \xrightarrow{\rho_{E,k}} E[f^{-1}]
\xrightarrow{} 0
$$is exact, where $\varepsilon_k^{-1}(P\otimes \delta
\otimes (f^{-k}e)) = (P\delta)\otimes (f^{-k}e) - P\otimes
\delta(f^{-k}e)$ (see \ref{nume:CECRSP}) and $\rho_{E,k}(P\otimes
(f^{-k}e)) = P(f^{-k}e)$.
\smallskip

\noindent From Remark \ref{nota:itera-BS}, b), we know that
$\rho_{E,k}$ is surjective if $k\geq -m_p$.\smallskip

\noindent The inclusion $\im \varepsilon_k^{-1} \subset \ker
\rho_{E,k}$ is clear. For the opposite inclusion, we know by
Proposition \ref{prop:kerrhos-kerrhok} that $\ker
\rho_{E,k}=\Phi_{E,k}^0(\ker \rho_{E,s})$ for any $k\geq -m_p$, and
so from Theorem \ref{teo:main-0} we obtain
\begin{eqnarray*}
& \ker \rho_{E,k} = \Phi_{E,k}^0\left(\ker \rho_{E,s}\right) =
\Phi_{E,k}^0\left(\im
\varepsilon_s^{-1}\right)= &\\
& = \im \left(\Phi_{E,k}^0\pcirc \varepsilon_s^{-1}\right) = \im
\left(\varepsilon_k^{-1}\pcirc \Phi_{E,k}^1\right) \subset \im
\varepsilon_k^{-1}. &
\end{eqnarray*}
\end{prueba}

\begin{corolario} \label{cor:main} Let us suppose that $D$ is of
linear jacobian type and let $\LL$ be the local system of horizontal
sections of $\E$ on $U=X-D$. Let $m_p$ be the smallest integer root
of $b_{\E,p}(s)$. Then, there is an open neighborhood $V$ of $p$
such that the restriction to $V$ of the canonical morphism
$$\Omega_X^{\bullet}(\log D)(\E(kD))\xrightarrow{} R j_* \LL$$
is an isomorphism in the derived category for $k\geq -m_p$.
\end{corolario}

\begin{prueba} It is a consequence of Theorems \ref{teo:crit-LCT}
and \ref{teo:main}.
\end{prueba}
\medskip

The above corollary answers a questions raised in
\cite{calde_nar_fourier}, Ex.~5.3.

\begin{corolario} \label{cor:main-bis} With the same hypothesis as Corollary \ref{cor:main},
let $m_p^*$ be the smallest integer root of $b_{\E^*,p}(s)$. Then,
there is an open neighborhood $V$ of $p$ such that the restriction
to $V$ of the canonical morphism
$$j_!\LL \xrightarrow{} \Omega_X^{\bullet}(\log D)(\E(rD))$$
is an isomorphism in the derived category for $r < m_p^*$.
\end{corolario}

\begin{prueba} It is a consequence of Corollary (3.1.8) in
\cite{calde_nar_fourier} and Corollary \ref{cor:main} applied to the
dual connection $\E^*$.
\end{prueba}

\begin{nota}
In Theorem \ref{teo:main} we obtain a global isomorphism ($V=X$) if
$m:=\inf_{p\in D} m_p > -\infty$ and $k\geq -m$. A similar remark
applies to Corollaries \ref{cor:main} and \ref{cor:main-bis}.
\end{nota}

\begin{nota} In the case $\E=\OX$, Theorem \ref{teo:main} would give
a proof of the LCT  provided that the Bernstein-Sato polynomial of a
reduced local equation of $D$ at any point has no integer roots less
than $-1$. This is the case when $D$ is locally quasi-homogeneous,
but we do not know any direct proof of this fact (see Remark
\ref{nota:coment-LCT}).
\end{nota}

\begin{nota} It would be interesting to have a proof of Theorem
\ref{teo:main} by using property 5) in Theorem \ref{teo:crit-LCT},
in a similar way as we did in Corollary \ref{cor:LCT} for the case
of locally quasi-homogeneous free divisors and $\E=\OX$.
\end{nota}

\begin{nota} Let $X$ be a smooth algebraic variety over $\C$, or over a field of
characteristic zero, and $D\subset X$ a hypersurface. The property
of being free, Koszul free, of linear jacobian type or of
differential linear type makes sense in the algebraic category, and
there is an algebraic version
of Theorem \ref{teo:main} whose proof
seems possible following the lines in this paper.
\end{nota}

\section{Applications to intersection $D$-modules}

In this section we assume that $D\subset X$ is a free divisor of
linear jacobian type, and $\E$ is an ILC with respect to $D$.
\medskip

Let $m_p$ (resp. $m_p^*$) be the smallest integer root of
$b_{\E,p}(s)$ (resp. of $b_{\E^*,p}(s)$), and let us assume that
$$ m:=\inf_{p\in D} m_p > -\infty\quad \text{and}\quad
m^*:=\inf_{p\in D} m_p^* > -\infty.  $$

Let $\LL$ be the local system of the horizontal sections of $\E$ on
$U=X-D$. As we saw in section \ref{subsec:LCP}, the canonical
morphism $ \DR \E(\star D)  =\Omega_X^{\bullet}(\E(\star D))
\xrightarrow{} R j_* \LL$ is an isomorphism in the derived category.
On the other hand, since $D$ is Koszul,
$$ \DX\Lotimes_{\DX(\log D)} \E(kD) = \DX\otimes_{\DX(\log D)}
\E(kD),
$$
and from Theorem \ref{teo:main} we deduce that
\begin{equation} \label{eq:RjL}
 \DR \left(\DX \otimes_{\DX(\log D)} \E(kD)\right)  \simeq
\DR \E(\star D) \simeq \Omega_X^{\bullet}(\E(\star D)) \simeq R j_*
\LL
\end{equation}
for  $k\geq -m$.
\smallskip

 Let us consider now the dual local system
$\LL^{\vee}$, which appears as the local system of the horizontal
sections of the dual ILC $\E^*$. Proceeding as above, we find that
\begin{equation} \label{eq:RjLdual}
 \DR \left(\DX \otimes_{\DX(\log D)} \E^*(k'D)\right)  \simeq
 \DR \E^*(\star D) \simeq \Omega_X^{\bullet}(\E^*(\star D)) \simeq R
j_* \LL^{\vee}
\end{equation}
for $k'\geq -m^*$.
\medskip

\noindent For $k+k'\geq 1$ let us denote by
$$ \varrho_{\E,k,k'}:\DX \otimes_{\DX(\log D)}
\E((1-k')D) \xrightarrow{} \DX \otimes_{\DX(\log D)} \E(kD)
$$
the $\DX$-linear map induced by the inclusion $\E((1-k')D) \subset
\E(kD)$, and by $\IC_X(\LL)$ the intersection complex of
Deligne-Goresky-MacPherson associated with $\LL$, which is described
as the intermediate direct image $j_{!*}\LL$, i.e. the image of
$j_!\LL \to R j_*\LL$ in the category of perverse sheaves (cf.
\cite{bbd_83}, Def. 1.4.22).

\begin{steorema} \label{teo:2} Under the above conditions,
we have a canonical isomorphism in the category of perverse sheaves
on $X$
$$ \IC_X(\LL) \simeq \DR \left( \im \varrho_{\E,k,k'}\right),$$
for $k\geq -m,k'\geq -m^*$ and $k+k'\geq 1$. In other words, the
``intersection $\DX$-module" corresponding to $\IC_X(\LL)$ by the
Riemann-Hilbert correspondence of Mebkhout-Kashiwara
\cite{kas_RH,meb_I_84,meb_II_84} is $\im \varrho_{\E,k,k'}$, for
$k,k'\gg 0$.
\end{steorema}

\begin{prueba}
By using our duality results (\cite{calde_nar_fourier}, \S 3) and
the Local Duality Theorem for holonomic $\DX$-modules
(\cite{meb_formalisme}, ch. I, Th. (4.3.1); see also
\cite{nar-ldt}), we obtain
\begin{eqnarray*}
& \DR \left(\DX \otimes_{\DX(\log D)} \E((1-k')D)\right) \simeq \DR
 \left(\DX \otimes_{\DX(\log D)} \left(\E^*(k'D)\right)^*(D) \right) \simeq &\\
& \DR \left( \Dual_{\DX} \left(\DX \otimes_{\DX(\log D)} \E^*(k'D)
\right) \right) \simeq \left[\DR \left( \DX \otimes_{\DX(\log D)}
\E^*(k'D) \right)\right]^{\vee}
\stackrel{\text{(\ref{eq:RjLdual})}}{\simeq}&\\ & \left[R j_*
\LL^{\vee}\right]^{\vee} \simeq j_! \LL.&
\end{eqnarray*}

On the other hand, the canonical morphism $j_! \LL \to R j_* \LL$
corresponds, through the de Rham functor, to the $\DX$-linear
morphism $\varrho_{\E,k,k'}$, and the theorem is a consequence of the
Riemann-Hilbert correspondence which says that the de Rham functor
establishes an equivalence of abelian categories between the category
of regular holonomic $\DX$-modules and the category of perverse
sheaves on $X$.
\end{prueba}
\medskip

In \cite{calde_nar_MEGA} we use Theorem \ref{teo:2} to perform
explicit computations in the case of locally quasi-homogeneous plane
curves.

\bigskip

{\small \noindent Departamento de \'{A}lgebra,
 Facultad de  Matem\'{a}ticas, Universidad de Sevilla, P.O. Box 1160, 41080
 Sevilla, Spain}. \\
{\small {\it E-mail}:  $\{$calderon,narvaez$\}$@algebra.us.es
 }

\end{document}